\documentclass[reqno,letterpaper,11pt]{amsart}
\pdfoutput=1
\usepackage{amssymb}
\usepackage{graphicx}
\usepackage{amsfonts}
\usepackage[margin=0.93in, letterpaper ]{geometry}
\usepackage{lineno}

\newtheorem{theorem}{Theorem}
\theoremstyle{plain}

\newtheorem{claim}[theorem]{Claim}

\newtheorem{conjecture}[theorem]{Conjecture}
\newtheorem{corollary}[theorem]{Corollary}

\newtheorem{definition}[theorem]{Definition}

\newtheorem{lemma}[theorem]{Lemma}

\newtheorem{proposition}[theorem]{Proposition}

\newtheorem{remark}[theorem]{Remark}

\numberwithin{equation}{section}
\numberwithin{theorem}{subsection}
\newcommand{\legendre}[2]{\genfrac{(}{)}{}{}{#1}{#2}}

\begin{document}
\title{\textsf{Small Representations of Finite Classical Groups\smallskip }}
\date{\textit{Submitted July 13, 2016. }}
\author{\textit{Shamgar Gurevich}}
\address{\textit{Department of Mathematics, University of Wisconsin,
Madison, WI 53706, USA.}}
\email{shamgar@math.wisc.edu}
\author{\textit{Roger Howe}}
\address{\textit{Department of Mathematics, Yale University, New Haven, CT
06520, USA.}}
\email{roger.howe@yale.edu}

\begin{abstract}
Finite group theorists have established many formulas that express
interesting properties of a finite group in terms of sums of characters of
the group. An obstacle to applying these formulas is lack of control over
the dimensions of representations of the group. In particular, the
representations of small dimensions tend to contribute the largest terms to
these sums, so a systematic knowledge of these small representations could
lead to proofs of important conjectures which are currently out of reach.
Despite the classification by Lusztig of the irreducible representations of
finite groups of Lie type, it seems that this aspect remains obscure. In
this note we develop a language which seems to be adequate for the
description of the\textbf{\ "small"} representations of finite classical
groups and puts in the forefront the notion of \textbf{rank}\textit{\ }of a
representation. We describe a method, the \textbf{"eta correspondence", }to
construct small representations, and we conjecture that our construction is
exhaustive. We also give a strong estimate on the dimension of small
representations in terms of their rank. For the sake of clarity, in this
note we describe in detail only the case of the finite symplectic groups.
\end{abstract}

\maketitle

\section{\textbf{Introduction}}

Finite group theorists have established formulas that enable expression of
interesting properties of a group $G$ in terms of quantitative statements on
sums of values of its characters. There are many examples \cite{Di, DS, LS,
LOST, M, Sa, Sh}. We describe a representative one. Consider the commutator
map%
\begin{equation}
\left[ ,\right] :G\times G\rightarrow G;\text{ \ \ }[x,y]=xyx^{-1}y^{-1},
\label{CM}
\end{equation}%
and for $g\in G$ denote by $\left[ ,\right] _{g}$ the set $\left[ ,\right]
_{g}=\{(x,y)\in G\times G;$ $[x,y]=g\}.$ In \cite{O} Ore conjectured that
for a finite non-commutative simple group $G$ the map (\ref{CM}) is onto,
i.e., $\#\left[ ,\right] _{g}\neq 0,$ for every $g\in G.$ The quantity $\#%
\left[ ,\right] _{g}$ is a class function on $G$ and Frobenius developed the
formula for its expansion as a linear combination of irreducible characters.
Frobenius' formula is 
\begin{equation}
\frac{\#[,]_{g}}{\#G}=1+\sum_{1\neq \rho \in Irr(G)}\frac{\chi _{\rho }(g)}{%
\dim (\rho )},  \label{O-Sum}
\end{equation}%
where for $\rho $ in the set $Irr(G)$ of isomorphism classes of irreducible
representations---aka irreps---of $G,$ we use the symbol $\chi _{\rho }$ to
denote its character. Estimating the sum in the right-hand side of (\ref%
{O-Sum}) for certain classes of elements in several important finite
classical groups was a major technical step in the recent proof \cite{LOST,
M} of the Ore conjecture. Given the Ore conjecture thus, the following
question naturally arises:\medskip

\textbf{Question. }What is the distribution of the commutator map (\ref{CM}%
)?\medskip

In \cite{Sh} Shalev conjectured that for a finite non-commutative simple
group $G$ the distribution of (\ref{CM}) is approximately uniform, i.e., 
\begin{equation}
\sum_{1\neq \rho \in Irr(G)}\frac{\chi _{\rho }(g)}{\dim (\rho )}=o(1),\text{
\ }g\neq 1,  \label{C-S}
\end{equation}%
in a well defined quantitative sense (e.g., as $q \rightarrow \infty$  for a finite non-commutative 
simple group of Lie type $ G=G (\mathbb{F}_q$)). This conjecture is wide open \cite{Sh}. It can be proven for the
finite symplectic group $Sp_{2}(\mathbb{F}_{q})$\footnote{%
For the rest of this note, $q$ is a power of an odd prime number $p$.}
invoking its explicit character table, and probably also for $Sp_{4}(\mathbb{%
F}_{q})$ \cite{Sr-1}. As was noted by Shalev in \cite{Sh}, one can verify
the uniformity conjecture for elements of $G$ with small centralizers, using
Schur's orthogonality relations for characters \cite{LOST, M, Sh}. 
\begin{figure}[ptb]\centering
\includegraphics[
natheight=1.5558in, natwidth=9.9998in, height=0.8051in, width=5.028in]
{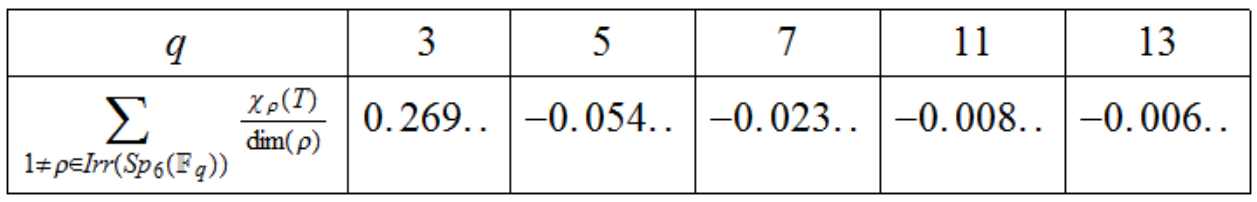}
\caption{Ore sum for the transvection $T$ (\protect\ref{Tr}) in $G=Sp_{6}(%
\mathbb{F}_{q})$ for various $q$'s.}\label{Ore}
\end{figure}
However, relatively little seems to be known about Shalev's conjecture in
the case of elements with relatively large centralizers---see Figure \ref%
{Ore} for numerical\footnote{%
The numerical data in this note was generated with J. Cannon (Sydney) and S.
Goldstein (Madison) using Magma.} illustration in the case of $G=Sp_{2n}(%
\mathbb{F}_{q})$ and the transvection element $T$ in $G$ which is given by 
\begin{equation}
T=%
\begin{pmatrix}
I & E \\ 
0 & I%
\end{pmatrix}%
,\text{ \ }E_{i,j}=\left\{ 
\begin{array}{c}
1,\text{ \ \ \ \ \ }i=j=1;\text{ \ \ \ \ \ \ } \\ 
0,\text{ \ \ other }1\leq i,j\leq n.%
\end{array}%
\right.  \label{Tr}
\end{equation}%
To suggest a possible approach for the resolution of the uniformity
conjecture, let us reinterpret (\ref{C-S}) as a statement about extensive
cancellation between the terms%
\begin{equation}
\frac{\chi _{\rho }(g)}{\dim (\rho )},\text{ }\rho \in Irr(G),  \label{C-R}
\end{equation}%
which are called \textit{character ratios. 
\begin{figure}[h]\centering
\includegraphics[
natheight=7.4996in, natwidth=9.9998in, height=4.2281in, width=5.6273in]
{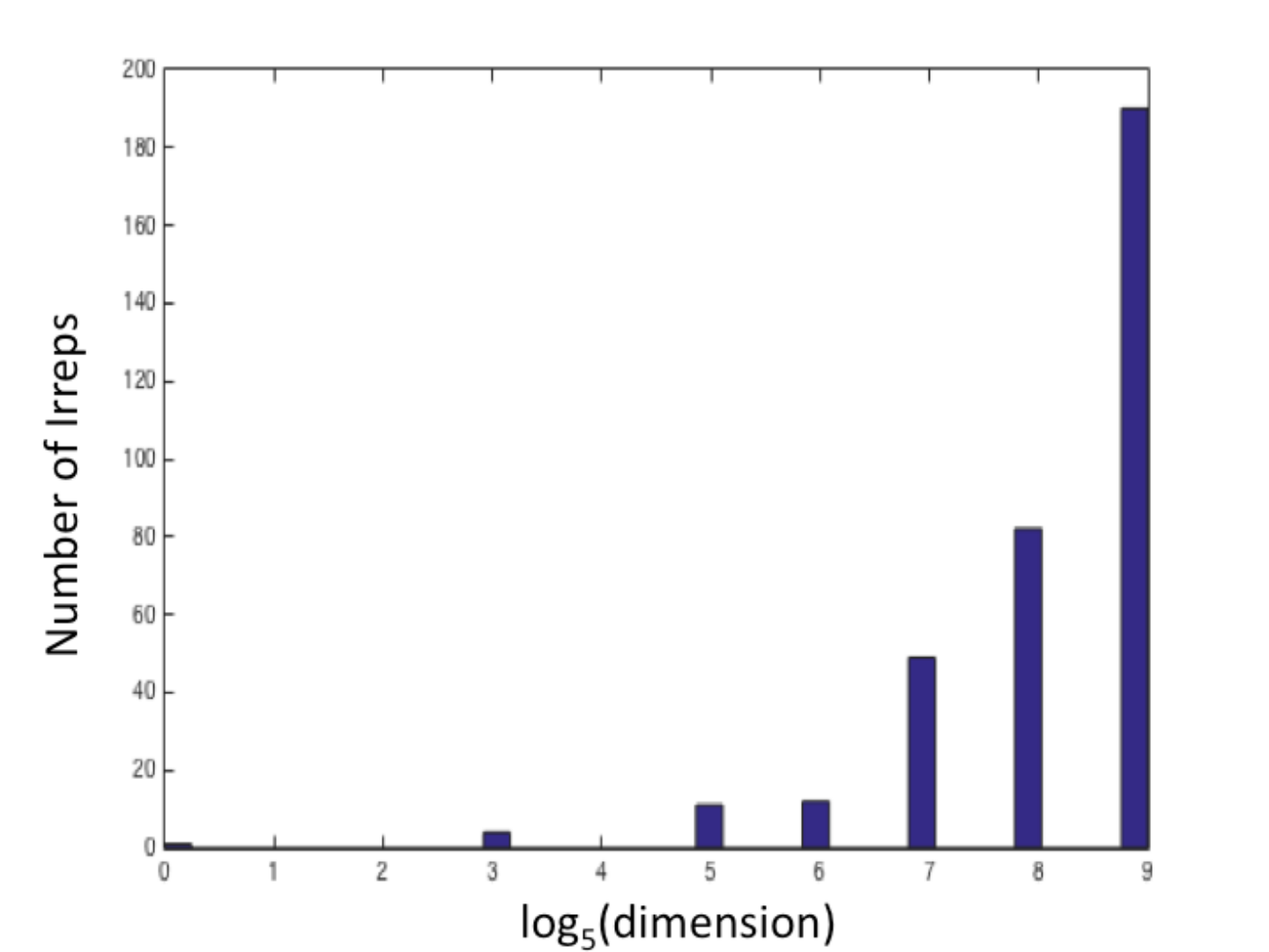}%
\caption{Partition of $Irr(Sp_{6}(\mathbb{F}_{5}))$ according to nearest
integer to $\log _{5}(\dim (\protect\rho )).$}\label{dim-irrsp6_5}%
\end{figure}%
}The dimensions of the irreducible representations of a finite group $G$
tend to come in certain layers according to order of magnitude. For
example---see Figure \ref{dim-irrsp6_5} for illustration---it is known \cite%
{DL, Lu1} that the dimensions of the irreducible representations of $%
G=Sp_{2n}(\mathbb{F}_{q})$ are given by some "universal" set of polynomials
in $q.$ In this case the degrees of these polynomials give a natural
partition of $Irr(Sp_{2n}(\mathbb{F}_{q}))$ according to order of magnitude
of dimensions$.$ Since the dimension of the representation of a group $G$ is
what appears in the denominator of (\ref{C-R}), it seems reasonable to
expect that\smallskip\ in (\ref{C-S})\smallskip

\textbf{(A) }\textit{Character ratios of lower dimensional representations
tend to contribute larger terms}.\smallskip

\textbf{(B) }\textit{The partial sums over low dimensional representations
of "similar" size exhibit cancellations.\medskip }

A significant amount of numerical data collected recently with Cannon and
Goldstein supports assertions (A) and (B). 
\begin{figure}[h]\centering
\includegraphics[
natheight=7.4996in, natwidth=9.9998in, height=4.9779in, width=6.6288in]
{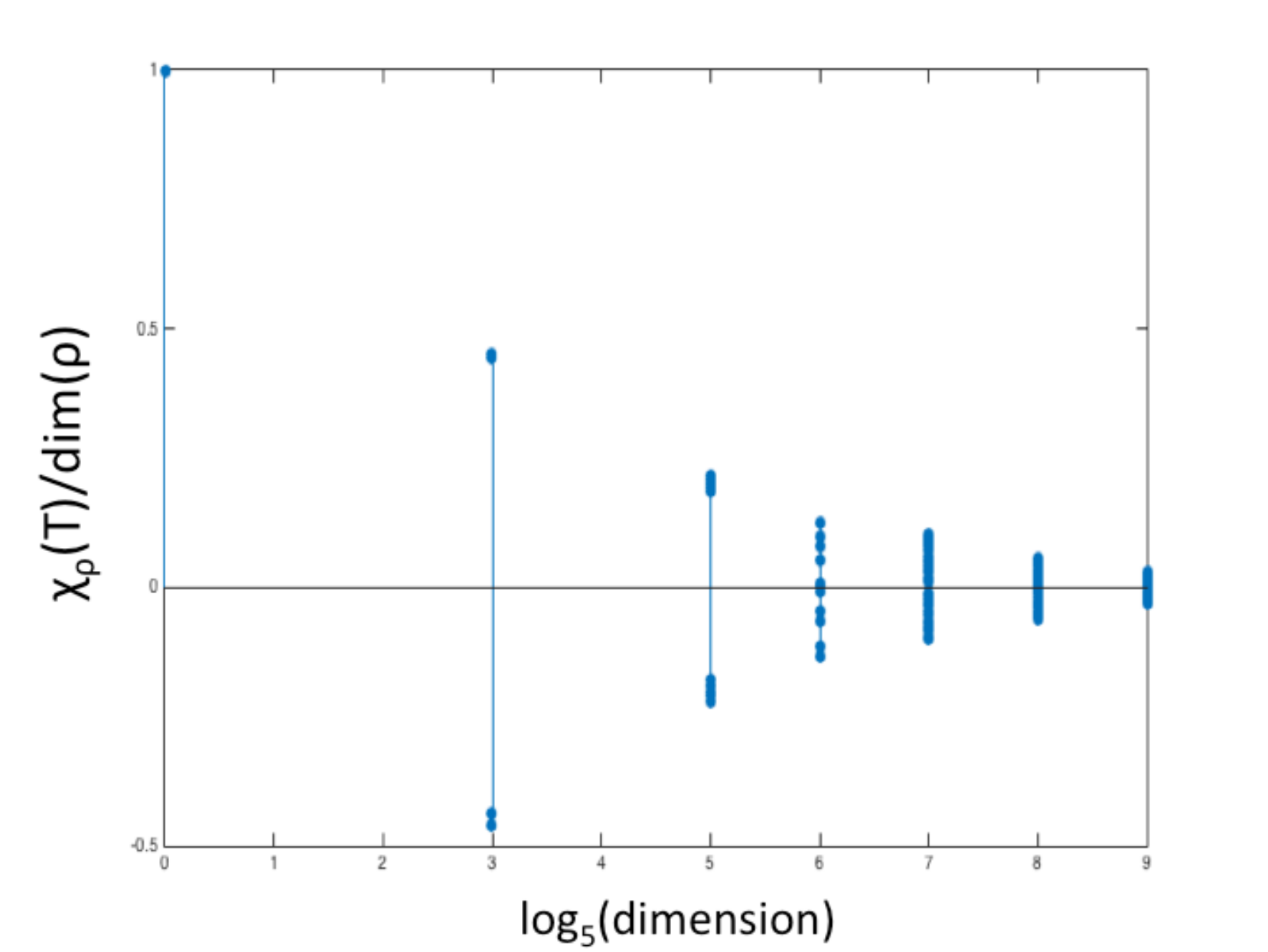}%
\caption{Character ratios at $T$ (\protect\ref{Tr})$\ $vs. nearest integer
to $\log _{5}(\dim (\protect\rho ))$ for $Irr(Sp_{6}(\mathbb{F}_{5}))$.}%
\label{CR-T-Sp6_5}%
\end{figure}%
For example, in Figure \ref{CR-T-Sp6_5} we plot the numerical values of the
character ratios of the irreducible representations of $G=Sp_{6}(\mathbb{F}%
_{5}),$ evaluated at the transvection $T$ (\ref{Tr}). More precisely, for
each $\rho \in Irr(Sp_{6}(\mathbb{F}_{5}))$ we marked by a circle the point%
\footnote{%
We denote by $\left\lfloor x\right\rceil $ the nearest integer value to a
real number $x.$} $(\left\lfloor \log _{5}(\dim (\rho )\right\rceil ,\chi
_{\rho }(T)/\dim (\rho ))$ and find that the overall picture is in agreement
with (A) and (B). Moreover---see Figures \ref{dim-irrsp6_5} and \ref%
{CR-T-Sp6_5} for illustration---the numerics shows that, although the
majority of representations are "large", their character ratios tend to be
so small that adding all of them contributes very little to the entire Ore
sum (\ref{C-S}). 

The above example illustrates that a possible obstacle to
getting group theoretical properties by summing over characters, as in
Formula (\ref{O-Sum}), is lack of control over the representations with
relatively small dimensions. In particular, it seems that a systematic
knowledge on the \textbf{"small" }representations of finite classical groups
could lead to proofs of some important open conjectures, which are currently
out of reach. However, relatively little seems to be known about these small
representations \cite{LS, M, Sh, TZ}.

In this note we develop a language suggesting that the small representations
of the finite classical groups can be systematically described by studying
their restrictions to unipotent subgroups, and especially, using the notion
of \textbf{rank} of a representation \cite{H3, L, S}. In addition, we
develop a new method, called the "\textbf{eta correspondence"}\textit{, }to
construct small representations. We conjecture that our construction is
exhaustive. Finally, we use our construction to give a strong estimate on
the dimension of the small representations in terms of their rank. For the
sake of clarity of exposition we treat in this note only the case of the
finite symplectic groups $Sp_{2n}(\mathbb{F}_{q})$.\medskip

\textbf{Acknowledgement. }We would like to thank John Cannon and Steve
Goldstein for the help in generating the numerical data appearing in this
note. We are grateful to Bob Guralnick for important conversations. Thanks
go also to the following institutions: MPI Bonn, Texas A\&M, UW-Madison,
Weizmann Institute, and Yale, where the work on this note was carried out.
Finally, S.G. thanks the UW-Madison for choosing him as a Vilas Associate
for the period of working on this project.

\section{\textbf{Notion of Rank of Representation}}

\begin{figure}[h]\centering
\includegraphics[
natheight=0.486in, natwidth=8.8332in, height=0.3589in, width=6.0347in]
{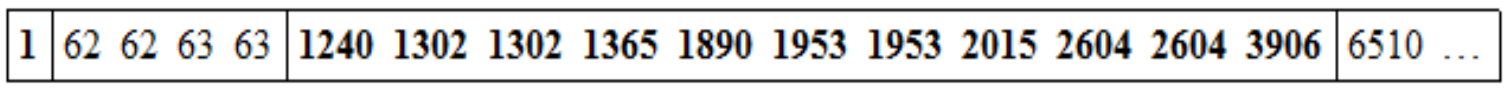}%
\caption{Dimensions of Irreps of $Sp_{6}(\mathbb{F}_{5}).$}%
\label{Numerical-Dim-IrrSp6_5}%
\end{figure}

Let us start with the numerical example of the dimensions of the irreducible
representations of the group $Sp_{6}(\mathbb{F}_{5}).$ The beginning of the
list appears in Figure \ref{Numerical-Dim-IrrSp6_5}. 
These numbers---see also Figure \ref{dim-irrsp6_5}---reveal the story of the
hierarchy in the world of representations of finite classical groups. A lot
of useful information is available on the "minimal" representations of these
groups, i.e., the ones of lowest dimensions \cite{LS, M, Sh, TZ}. In the
case of $Sp_{2n}(\mathbb{F}_{q})$ these are the $4$ components of the
oscillator (aka Weil) representations\textit{\ }\cite{G, GH, H1, H2, W}, $2$
of dimension $(q^{n}-1)/2$ and $2$ of dimension $(q^{n}+1)/2,$ which in
Figure \ref{Numerical-Dim-IrrSp6_5} are the ones of dimensions $62,62,63,63$%
. In addition, a lot is known about the "big" representations of the finite
classical groups, i.e., those of considerably large dimension (See \cite{C,
DL, DM, LMT, LS, Lu1, Lu2, Lu3, Sh, Sr1} and references therein). We will
not attempt to define the "big" representations at this stage, but in Figure %
\ref{Numerical-Dim-IrrSp6_5} the ones of dimension $6510$ and above fall in
that category. However, relatively little seems to be known about the
representations of the classical groups which are in the range between
"minimal" and "big" \cite{LS, M, Sh, TZ}. \ In Figure \ref%
{Numerical-Dim-IrrSp6_5} those form the layer of $11$ representations of
dimensions between $1240$ and $3906.$

In this section we introduce a language that will enable us to define the
"small" representations of finite classical groups. This language will
extend well beyond the notion of minimal representations and will induce a
partition of the set of isomorphism classes of irreducible representations
which is closely related to the hierarchy afforded by dimension. In
particular, this language gives an explicit organization of all the
representations in Figure \ref{Numerical-Dim-IrrSp6_5}, and explains why
this list is, in a suitable sense, complete. The key idea we will use is
that of the \textit{rank} of a representation. This notion was developed in
the 1980s by Howe, in the context of unitary representations of classical
groups over local fields \cite{H4}, but it has not been applied to finite
groups. For the sake of clarity of exposition, in this note we give the
definition of rank only in the symplectic case, leaving the more general
treatment to future publication. We start by discussing necessary
ingredients from the structure theory of $Sp_{2n}(\mathbb{F}_{q})$.

\subsection{\textbf{The Siegel Unipotent Radical}}

\ Let $(V,\left\langle ,\right\rangle )$ be a $2n$-dimensional symplectic
vector space over the finite field $\mathbb{F}_{q}.$ In order to simplify
certain formulas, let us assume that 
\begin{equation}
V=X\oplus Y,  \label{LD}
\end{equation}%
where $X$ and $Y$ are vector spaces dual to each other with pairing $\bullet 
$, and that the symplectic form $\left\langle ,\right\rangle $ is the
natural one which is defined by that pairing, i.e.,

\begin{equation}
\left\langle 
\begin{pmatrix}
x_{1} \\ 
y_{1}%
\end{pmatrix}%
,%
\begin{pmatrix}
x_{2} \\ 
y_{2}%
\end{pmatrix}%
\right\rangle =x_{1}\bullet y_{2}-x_{2}\bullet y_{1}.  \label{SympForm}
\end{equation}%
Note that $X$ and $Y$ are maximal isotropic---aka \textit{Lagrangian}%
---subspaces of $V.$ Consider the symplectic group $Sp=Sp(V)$ of elements of 
$GL(V)$ which preserve the form $\left\langle ,\right\rangle .$ Denote by $%
P=P_{X}$ the subgroup of all elements in $Sp$ that preserve $X$\textit{. }\
The group $P$ is called the \textit{Siegel parabolic} \cite{Si} and can be
described explicitly in terms of the decomposition (\ref{LD}) 
\begin{equation*}
P=\left\{ 
\begin{pmatrix}
I & A \\ 
0 & I%
\end{pmatrix}%
\cdot 
\begin{pmatrix}
C & 0 \\ 
0 & ^{t}C^{-1}%
\end{pmatrix}%
;\text{ }A:Y\rightarrow X\text{ \ symmetric, }C\in GL(X)\right\} ,
\end{equation*}%
where $^{t}C^{-1}\in GL(Y)$ is the inverse of the transpose of $C.$ In
particular, $P$ has the form\textit{\ }of a semi-direct product, known also
as its Levi decomposition \cite{DM}, 
\begin{equation}
P\simeq N\rtimes GL(X),  \label{NGL}
\end{equation}%
where $N=N_{X}$, called the \textit{unipotent radical} of $P$, \ is the
normal subgroup 
\begin{equation*}
N=\left\{ 
\begin{pmatrix}
I & A \\ 
0 & I%
\end{pmatrix}%
;\text{ }A:Y\rightarrow X\text{ \ symmetric}\right\} .
\end{equation*}%
The group $N$ is abelian and we have a tautological $GL(X)$-equivariant
isomorphism 
\begin{equation}
N\widetilde{\longrightarrow }Sym^{2}(X),  \label{SymX}
\end{equation}%
where $Sym^{2}(X)$ denotes the space of symmetric bilinear forms on $%
Y=X^{\ast },$ and the $GL(X)$ action on $Sym^{2}(X)$ is the standard one. In
addition, if we fix a non-trivial additive character $\psi $ of $\mathbb{F}%
_{q}$ we obtain a $GL(X)$-equivariant isomorphism%
\begin{equation}
\left\{ 
\begin{array}{c}
Sym^{2}(Y)\widetilde{\longrightarrow }\widehat{N}, \\ 
B\mapsto \psi _{B},%
\end{array}%
\right.  \label{SymXdual}
\end{equation}%
where $Sym^{2}(Y)$ denotes the space of symmetric bilinear forms on $%
X=Y^{\ast },$ the $GL(X)$ action on $Sym^{2}(Y)$ is the standard one, the
symbol $\widehat{N}$ stands for the Pontryagin dual (group of characters) of 
$N,$ and 
\begin{equation}
\psi _{B}(A)=\psi (Tr(BA)),  \label{Char-SymForms}
\end{equation}%
\ for every $A\in Sym^{2}(X),$ where $Tr(BA)$ indicates the trace of the
composite operator $Y\overset{A}{\rightarrow }X\overset{B}{\rightarrow }Y.$

\subsection{\textbf{The }$N$\textbf{-spectrum of a Representation}}

Now, take a representation $\rho $ of $Sp$ and look at the
restriction to $N$. It decomposes \cite{T} as a sum of characters with
certain multiplicities%
\begin{equation}
\rho _{|N}=\sum_{B\in Sym^{2}(Y)}m_{B}\psi _{B}.  \label{FT1}
\end{equation}%
The function $m$ and its support will be called, respectively, the 
\underline{$N$-spectrum} of $\rho $, and the \underline{$N$-support} of $%
\rho ,$ and will be denoted by $Spec_{N}(\rho ),$ and $Supp_{N}(\rho ),$
respectively.

We would like to organize the decomposition (\ref{FT1}) in a more meaningful
way. Note that the restriction to $N$ of a representation $\rho $ of $Sp$
can be thought of as the restriction to $N$ of the restriction of $\rho $ to 
$P$. Using (\ref{NGL}), this implies \cite{Ma1}:

\begin{proposition}
The $N$-spectrum of a representation $\rho $ of $Sp$ is $GL(X)$ invariant.
That is, $m_{B}=m_{B^{\prime }}$ if $B$ and $B^{\prime }$\ define equivalent
symmetric bilinear forms on $X$.
\end{proposition}

The first major invariant of a symmetric bilinear form is its rank. It is
well known \cite{La} that, over finite fields of odd characteristics, there
are just two isomorphism classes of symmetric bilinear forms of a given rank 
$r.$ They are classified by their discriminant \cite{La}, which is an
element in $\mathbb{F}_{q}^{\ast }/\mathbb{F}_{q}^{\ast 2}$. We denote by $%
\mathcal{O}_{r+}$ and $\mathcal{O}_{r-}$, the two classes of symmetric
bilinear forms, these whose discriminant is the coset of squares, and these
whose coset consists of non-squares, respectively; or we will denote the
pair of them, or whichever one is relevant in a given context as $\mathcal{O}%
_{r\pm }.$ If $B$ is a form of rank $r$, we will also say that the
associated character $\psi _{B}$ has rank $r$. We may also refer to the
character as being of type $+$ or type $-$, according to the type of $B$.
With this notation, we can reorganize the expansion (\ref{FT1}) of $\rho
_{|N}.$ Namely, we split the sum according to the ranks of characters, and
within each rank we split the sum into two partial sub-sums according to the
two isomorphism classes of the associated forms:%
\begin{equation}
\rho _{|N}=\sum_{r}\sum_{\pm }m_{r\pm }\sum_{B\in \mathcal{O}_{r\pm }}\psi
_{B}.  \label{FT2}
\end{equation}%
Note that, Formula (\ref{FT2}) implies, by evaluation at the identity of $N$%
, that the dimension of $\rho $ must be%
\begin{equation}
\dim (\rho )=\sum_{r}\sum_{\pm }m_{r\pm }\cdot \#\mathcal{O}_{r\pm },
\label{DimForm}
\end{equation}%
i.e., a weighted sum of the cardinalities $\#\mathcal{O}_{r\pm }$ of the
isomorphism classes of symmetric bilinear forms. It is easy to write
formulas for these cardinalities \cite{A}. We have%
\begin{equation}
\#\mathcal{O}_{r\pm }=\#Gr_{n,r}\cdot \frac{\#GL_{r}}{\#O_{r\pm }},
\label{Card1}
\end{equation}%
where $Gr_{n,r}$ denotes the Grassmannian of $r$-dimensional subspaces of $%
\mathbb{F}_{q}^{n}$, the symbol $GL_{r}$ stands for the group of
automorphisms of $\mathbb{F}_{q}^{r},$ and $O_{r\pm }$ is the isometry group
of a non-degenerate form of type $\pm $ on $\mathbb{F}_{q}^{r},$ i.e., it is 
$O_{r+}$ in case of a form from $\mathcal{O}_{r+}$ and likewise with $+$ and 
$-$ switched. In particular, using standard formulas \cite{A} for $%
\#Gr_{n,r},$ $\#GL_{r},$ and $\#O_{r\pm },$ we obtain 
\begin{equation}
\#\mathcal{O}_{r\pm }\approx \frac{1}{2}q^{r(n-\frac{r-1}{2})}\text{ .}
\label{Card2}
\end{equation}

\subsection{\textbf{Smallest Possible Irreducible Representation}}

From (\ref{Card2}) we get, in particular, that the smallest non-trivial
orbits are those of rank one forms. Using (\ref{Card1}) we see that these
have size $\#\mathcal{O}_{1\pm }=(q^{n}-1)/2.$ It follows from this that the
smallest possible dimension of a non-trivial irreducible representation $%
\rho $ of $Sp$ should satisfy%
\begin{equation}
\dim (\rho )\geq \frac{q^{n}-1}{2}.  \label{LowestDim}
\end{equation}%
Indeed, we have the following lemma:

\begin{lemma}
\label{Trivial}The only irreducible representation of $Sp$ with $N$-spectrum
concentrated at zero is the trivial one.
\end{lemma}

The proof of Lemma \ref{Trivial} is easy, but to avoid interrupting this
discussion, we defer it to Appendix \ref{P-of-Lemma-Trivial}.\medskip

A representation whose dimension attaining the lower bound (\ref{LowestDim})
would contain each rank one character of one type, and nothing else. Since $%
N $ is such a small subgroup of $Sp$, it is unclear whether to expect such a
representation to exist. In particular, it would be irreducible already on
the Siegel parabolic $P$, and it would be the smallest possible faithful
representation of $P$. It turns out, however, that it does exist; in fact,
there are two \cite{G, H1, H2, I, W}.

\begin{proposition}
\label{Smallest}There are two irreducible representations of $Sp$ of
dimension $\frac{q^{n}-1}{2}$, one containing either one of the two rank one 
$GL(X)$ orbits in $\widehat{N}.$
\end{proposition}

What is the next largest possible dimension? Well, one more - the $N$%
-spectrum could include a rank one orbit, and a trivial representation. It
turns out that these also exist \cite{G, H1, H2, I, W}.

\begin{proposition}
\label{Next}There are two irreducible representations of $Sp$ of dimension $%
\frac{q^{n}-1}{2}+1=\frac{q^{n}+1}{2}$ one whose $N$-spectrum contains one
of the rank one orbits in $\widehat{N}.$
\end{proposition}

For a proof of Propositions \ref{Smallest} and \ref{Next}, see Section \ref%
{Smallest Rep}.

\subsection{\textbf{Definition of Rank of Representation}}

The existence of the above smallest possible representations, plus
considerations of tensor products, tell us that, for any orbit $\mathcal{O}%
_{k\pm }$ in $\widehat{N}$ there will be representations of $Sp$ whose $N$%
-spectrum contains the given orbit, together with orbits of smaller rank.
Since the size---see Formulas (\ref{Card1}) and (\ref{Card2})---of the
orbits $\mathcal{O}_{k\pm }$ is increasing rapidly with $k$, representations
whose $N$-spectrum is concentrated on orbits of smaller rank can be expected
to have smaller dimensions. This motivates us to introduce the following key
notion in our approach for small representations.

\begin{definition}[\textbf{Rank}]
Let $\rho $ be a representation of $Sp.$

\begin{enumerate}
\item We say that $\rho $ is of \underline{rank} $k,$ denoted $rk(\rho )=k,$
iff the restriction $\rho _{|N}$ contains characters of rank $k$, but of no
higher rank.

\item If $\rho $ is of rank $k$ and contains characters of type $\mathcal{O}%
_{k+},$ but not of type $\mathcal{O}_{k-}$, then we say that $\rho $ is of 
\underline{type} $\mathcal{O}_{k+}$; and likewise with $+$ and $-$ switched.
\end{enumerate}
\end{definition}

\begin{figure}[h]\centering
\includegraphics[
natheight=4.875in, natwidth=16.9166in, height=1.9787in, width=6.7948in]
{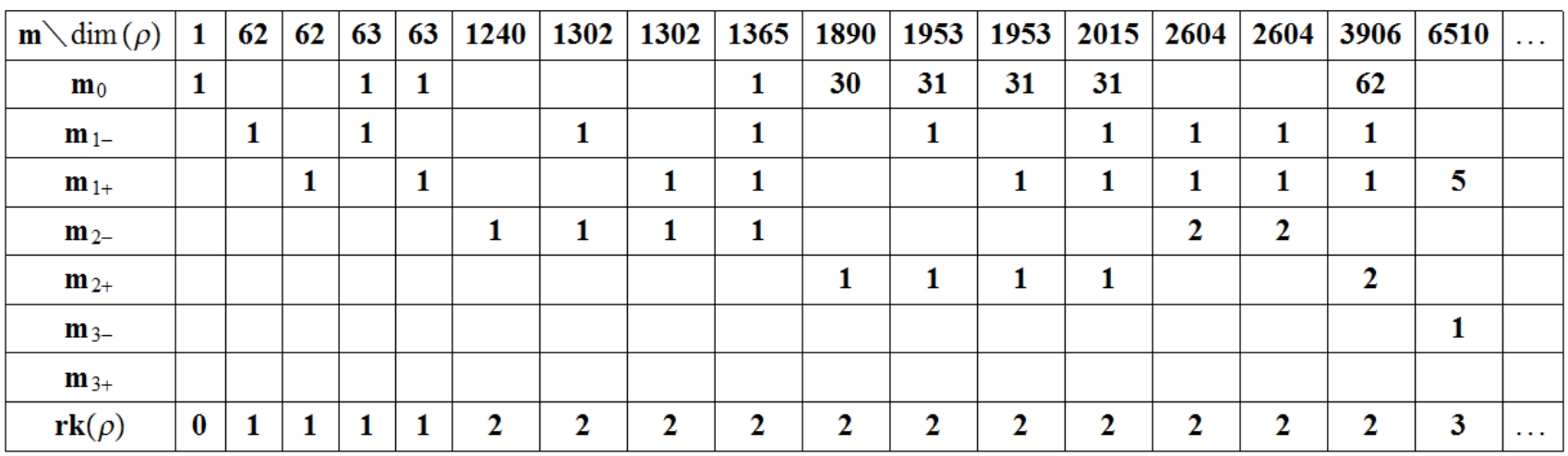}
\caption{Multiplicities and rank for irreps of $Sp_{6}(\mathbb{F}_{5}).$}%
\label{RankSp6_5}
\end{figure}

Let us convey some intuition for this notion using numerical data obtained
for the irreducible representations of the group $Sp_{6}(\mathbb{F}_{5})$%
---see Figure \ref{RankSp6_5}. The computations of the multiplicities and
rank in this case reveal a striking compatibility with the families of
representations appearing in the list of Figure \ref{Numerical-Dim-IrrSp6_5}%
. 
For example, it shows that the trivial representation is the one with rank $%
k=0;$ the $4$ components of the two oscillator representations are those of
rank $k=1$ and they split into $2$ of type $\mathcal{O}_{1+}$ and $2$ of
type $\mathcal{O}_{1-};$ the $11$ representations of dimensions between $%
1240 $ and $3906$ are the ones of rank $k=2$ and they split into $5$ of type 
$\mathcal{O}_{2+}$ and $6$ of type $\mathcal{O}_{2-}$; and above that the
"big" representations are those with rank $k=3.$

The main quest now is for a systematic construction of the "low rank"
irreducible representations. In the next section we take the first step
toward that goal by treating the smallest non-trivial representations of $Sp$
which are of rank $k=1$---see Propositions \ref{Smallest} and \ref{Next}.

\section{\textbf{The Heisenberg and Oscillator Representations}}

Where do the smallest representations of $Sp$ come from? A conceptual answer
to this question was given by Weil in \cite{W}. They can be found by
considering the Heisenberg\textit{\ }group.

\subsection{\textbf{The Heisenberg Group}}

The \textit{Heisenberg }group attached to $(V,\left\langle ,\right\rangle )$
is a two-step nilpotent group that can be realized by the set%
\begin{equation*}
H=V\times \mathbb{F}_{q},
\end{equation*}%
with the group law%
\begin{equation*}
\left( v,z\right) \cdot (v^{\prime },z^{\prime })=(v+v^{\prime },z+z^{\prime
}+\frac{1}{2}\left\langle v,v^{\prime }\right\rangle ).
\end{equation*}%
In particular, the center $Z$ of the Heisenberg group 
\begin{equation*}
Z=\{(0,z);\text{ }z\in \mathbb{F}_{q}\},
\end{equation*}%
is equal to its commutator subgroup. Moreover, the commutator operation in $%
H $ induces a skew-symmetric bilinear form on $H/Z\simeq V$ that coincides
with the original symplectic form.

The group $H$ is the analog over a finite field of the Lie group associated
with the Canonical Commutation Relations (CCR) of Werner Heisenberg, of
Uncertainty Principle fame.

\subsection{\textbf{Representations of the Heisenberg Group}}

We would like to describe the representation theory, i.e., the irreducible
representations, of the Heisenberg group. This theory is simultaneously
simple and deep, with fundamental connections to a wide range of areas in
mathematics and its applications. Take an irreducible representation $\pi $
of $H$. Then, by Schur's lemma, the center $Z$ will act by scalars%
\begin{equation*}
\pi (0,z)=\psi _{\pi }(z)I,\text{ \ }z\in Z,
\end{equation*}%
where $I$ is the the identity operator on the representation space of $\pi $%
, and $\psi _{\pi }\in \widehat{Z}$ is a character of $Z,$ called the 
\textit{central character }of $\pi $. If $\psi _{\pi }=1,$ then $\pi $
factors through $H/Z\simeq V$, which is abelian, so $\pi $ is itself a
character of $V.$ The case of non-trivial central character is described by
the following celebrated theorem \cite{Ma}:

\begin{theorem}[\textbf{Stone--von Neumann--Mackey}]
Up to equivalence, there is a unique irreducible representation $\pi _{\psi
} $ with given non-trivial central character $\psi $ in $\widehat{Z}%
\smallsetminus \{1\}.$
\end{theorem}

We will call the (isomorphism class of the) representation $\pi _{\psi }$
the \textit{Heisenberg representation }associated to the central character $%
\psi .$

\begin{remark}[\textbf{Realization}]
There are many ways to realize (i.e., to write explicit formulas for) $\pi
_{\psi }$ \cite{G, GH, GH2, H1, H2, I, W}. In particular, it can be
constructed as induced representation from any character extending $\psi $
to any maximal abelian subgroup of $H$ \cite{H3.5, Ma1}. To have a concrete
one, note that the inverse image in $H$ of any Lagrangian subspace of $V$
will be a maximal abelian subgroup for which we can naturally extend the
character $\psi $. For example, consider the Lagrangian $X\subset V$ and the
associated maximal abelian subgroup $\widetilde{X}$ with character $%
\widetilde{\psi }$ on it, given by 
\begin{equation*}
\widetilde{X}=X\times \mathbb{F}_{q},\text{ \ }\widetilde{\psi }(x,z)=\psi
(z).
\end{equation*}%
Then we have the explicit realization of $\pi _{\psi },$ given by the action
of $H,$ by right translations, on the space 
\begin{equation}
Ind_{\widetilde{X}}^{H}(\widetilde{\psi })=\{f:H\rightarrow 
\mathbb{C}
;\text{ \ }f(\widetilde{x}h)=\widetilde{\psi }(\widetilde{x})f(h)\text{, \ \ 
}\widetilde{x}\in \widetilde{X}\text{, }h\in H\}.  \label{Ind}
\end{equation}%
In particular, we have $\dim (\pi _{\psi })=q^{n}.$
\end{remark}

\subsection{\textbf{The Oscillator Representation}}

A compelling property of the Heisenberg group is that it has a large
automorphism group. In particular, the action of $Sp$ on $V$ lifts to an
action on $H$ by automorphisms leaving the center point-wise fixed. The
precise formula is $g(v,z)=(gv,z),$ \ $g\in Sp.$ It follows from the \textit{%
Stone--von Neumann--Mackey} theorem, that the induced action of $Sp$ on the
set $Irr(H)$ will leave fixed each isomorphism class $\pi _{\psi },$ $\psi
\in \widehat{Z}\smallsetminus \{1\}.$ This means that, if we fix a vector
space $\mathcal{H}_{\psi }$ realizing $\pi _{\psi },$ then for each $g$ in $%
Sp$ there is an operator $\omega _{\psi }(g)$ which acts on space $\mathcal{H%
}_{\psi }$ and satisfies the equation 
\begin{equation}
\omega _{\psi }(g)\pi _{\psi }(h)\omega _{\psi }(g)^{-1}=\pi _{\psi }(g(h)),
\label{Egorov}
\end{equation}%
which is also known as the exact \textit{Egorov identity} \cite{E} in the
mathematical physics literature. Note that, by Schur's lemma, the operator $%
\omega _{\psi }(g)$ is defined by (\ref{Egorov}) up to scalar multiples.
This implies that for any $g,g^{\prime }\in Sp$ we have $\omega _{\psi
}(g)\omega _{\psi }(g^{\prime })=c(g,g^{\prime })\omega _{\psi }(gg^{\prime
}),$ where $c(g,g^{\prime })$ is an appropriate complex number of absolute
value $1$. It is well known (see \cite{G, GH, GH2} for explicit formulas)
that over finite fields of odd characteristic this mapping can be lifted to
a genuine representation.

\begin{theorem}[\textbf{Oscillator Representation}]
\label{OR}There exists\footnote{%
The lift is unique except the case $n=2$ and $q=3,$ where still there is a
canonical lift \cite{GH, GH2}.} a representation 
\begin{equation*}
\omega _{\psi }:Sp\longrightarrow GL(\mathcal{H}),
\end{equation*}%
that satisfies the Egorov identity (\ref{Egorov}).
\end{theorem}

We will call $\omega _{\psi }$ the \textit{oscillator }representation. This
is a name that was given to this representation in \cite{H2} due to its
origin in physics \cite{Se, Sha}. Another popular name for $\omega _{\psi }$
is the \textit{Weil }representation, following the influential paper \cite{W}%
.

\begin{remark}[\textbf{Schr\"{o}dinger Model}]
\label{Schrodinger Model}We would like to have some useful formulas for the
representation $\omega _{\psi }$. Note that the space (\ref{Ind}) is
naturally identified with%
\begin{equation}
L^{2}(Y)\text{ \ - \ functions on }Y.  \label{Space}
\end{equation}%
On the space (\ref{Space}) we realize the representation $\omega _{\psi }$.
This realization is sometime called the \textit{Schr\"{o}dinger model. }In
particular, in that model for every $f\in L^{2}(Y)$ we have \cite{G, W}%
\bigskip

(A) $\left[ \omega _{\psi }%
\begin{pmatrix}
I & A \\ 
0 & I%
\end{pmatrix}%
f\right] (y)=\psi (\frac{1}{2}A(y,y))f(y),\medskip $

where $A:Y\rightarrow X$ is symmetric;\bigskip

(B) $\left[ \omega _{\psi }%
\begin{pmatrix}
0 & B \\ 
-B^{-1} & 0%
\end{pmatrix}%
f\right] (y)=\frac{1}{\gamma (B,\psi )}\sum\limits_{y^{\prime }\in Y}\psi
(B(y,y^{\prime }))f(y^{\prime }),\medskip $

where $B:Y\widetilde{\rightarrow }X$ is symmetric, and $\gamma (B,\psi
)=\sum\limits_{y\in Y}\psi (-\frac{1}{2}B(y,y))$ the quadratic Gauss
sum;\bigskip

(C) $\left[ \omega _{\psi }%
\begin{pmatrix}
^{t}C^{-1} & 0 \\ 
0 & C%
\end{pmatrix}%
f\right] (y)=\legendre{\det (C)}{q}f(C^{-1}y),\medskip $

where $C\in GL(Y),$ $^{t}C^{-1}\in GL(X)$ its transpose inverse, and $%
\legendre{\cdot}{q}$ is the Legendre symbol\footnote{%
For $x\in \mathbb{F}_{q}^{\ast }$ the Legendre symbol $\legendre{x}{q}=+1$
or $-1,$ according to $x$ being a square or not, respectively$.$}.
\end{remark}

It turns out that the isomorphism class of $\omega _{\psi }$ does change
when varying the central character $\psi $ in $\widehat{Z}\smallsetminus
\{1\}.$ However, this dependence is weak. The following result indicates
that there are only two possible oscillator representations. For a character 
$\psi $ in $\widehat{Z}\smallsetminus \{1\}$ denote by $\psi _{a},$ $a\in 
\mathbb{F}_{q}^{\ast },$ the character $\psi _{a}(0,z)=\psi (0,az).$

\begin{proposition}
\label{Oscillators}We have $\omega _{\psi }\simeq \omega _{\psi ^{\prime }}$
iff $\psi ^{\prime }=\psi _{s^{2}\text{ }}$ for some $s\in \mathbb{F}%
_{q}^{\ast }.$
\end{proposition}

For a proof of Proposition \ref{Oscillators}, see Appendix \ref%
{P-of-Proposition-Oscillators}.

\subsection{\textbf{The Smallest Possible Representations\label{Smallest Rep}%
}}

Using Formula (A) given in Remark \ref{Schrodinger Model}, it is easy to
determine the rank of the oscillator representation.

\begin{proposition}
Each representation $\omega _{\psi }$ is of rank $1$. One isomorphism class
is of type $\mathcal{O}_{1+}$ and the other is of type $\mathcal{O}_{1-}$.
\end{proposition}

In addition, the oscillator representations are slightly reducible. The
center $Z(Sp)=\{\pm I\}$ acts on the representation $\omega _{\psi }$---see
Remark \ref{Schrodinger Model} for the explicit action of $-I.$ We have the
direct sum decomposition 
\begin{equation}
\omega _{\psi }=\omega _{\psi ,1}\oplus \omega _{\psi ,sgn},
\label{Even-Odd}
\end{equation}%
with 
\begin{equation*}
\dim (\omega _{\psi ,1})=\left\{ 
\begin{array}{c}
\frac{q^{n}+1}{2}\text{ \ if \ }q\equiv 1\text{ }\text{mod}\text{ }4; \\ 
\frac{q^{n}-1}{2}\text{ \ if \ }q\equiv 3\text{ }\text{mod}\text{ }4;%
\end{array}%
\right. \text{ \ and\ \ }\dim (\omega _{\psi ,sgn})=\left\{ 
\begin{array}{c}
\frac{q^{n}-1}{2}\text{ \ if \ }q\equiv 1\text{ }\text{mod}\text{ }4; \\ 
\frac{q^{n}+1}{2}\text{ \ if \ }q\equiv 3\text{ }\text{mod}\text{ }4,%
\end{array}%
\right. \text{ }
\end{equation*}%
where $\omega _{\psi ,1}$ is the subspace of "even vectors", i.e., vectors
on which $Z(Sp)$ acts trivially, and $\omega _{\psi ,sgn}$ is the subspace
of "odd vectors", i.e., vectors on which $Z(Sp)$ acts via the sign
character. The above discussion also implies the following:

\begin{theorem}
The decomposition (\ref{Even-Odd}) is the decomposition of $\omega _{\psi }$
into irreducible representations.
\end{theorem}

To conclude, our study of the oscillator representation has established
Propositions \ref{Smallest} and \ref{Next}. More precisely, the
representations (\ref{Even-Odd}) have rank one, they are of type $\mathcal{O}%
_{1\pm }$, and have the required dimensions.

\section{\textbf{Construction of Rank }$k$\textbf{\ Representations}}

Where do higher rank representations of $Sp$ come from?\textbf{\ }This
section will include an answer to this question in the regime of "small"
representations. More, precisely we give here a systematic construction of
rank $k$ irreducible representations of $Sp$ in the so-called "stable range"%
\textit{\ }%
\begin{equation*}
k<n=\frac{\dim (V)}{2}.
\end{equation*}%
We will also refer to such representations as "small" or "low rank".

\subsection{The \textbf{Symplectic-Orthogonal Dual Pair}}

Let $U$ be a $k$-dimensional vector space over $\mathbb{F}_{q},$ and let $%
\beta $ be an inner product (i.e., a non-degenerate symmetric bilinear form)
on $U$. The pair $(U,\beta )$ is called a \textit{quadratic space. }We
denote by $O_{\beta }$ the isometry group of the form $\beta .$ Consider the
vector space $V\otimes U$---the tensor product of $V$ and $U$ \cite{AMa}. It
has a natural structure of a symplectic space, with the symplectic form
given by $\left\langle ,\right\rangle \otimes \beta .$ The groups $Sp=Sp(V)$
and $O_{\beta }$ act on $V\otimes U$ via their actions on the first and
second factors, respectively,%
\begin{equation*}
Sp\curvearrowright V\otimes U\curvearrowleft O_{\beta }.
\end{equation*}%
Both actions preserve the form $\left\langle ,\right\rangle \otimes \beta ,$
and moreover the action of $Sp$ commutes with that of $O_{\beta }$, and vice
versa. In particular, we have a map 
\begin{equation}
Sp\times O_{\beta }\longrightarrow Sp(V\otimes U),  \label{Embedding}
\end{equation}%
which embeds the two factors $Sp$ and $O_{\beta }$ in $Sp(V\otimes U),$ and
they form a pair of commuting subgroups. In fact, each is the full
centralizer of the other inside $Sp(V\otimes U).$ Thus, the pair $%
(Sp,O_{\beta })$ forms what has been called in \cite{H2} a \textit{dual pair}
of subgroups of $Sp(V\otimes U)$.

\subsection{The \textbf{Schr\"{o}dinger Model\label{Schrodinger-VU}}}

We write down a specific model for the oscillator representation\footnote{%
We suppress the dependence of $\omega _{V\otimes U}$ on the central
character, but we record which symplectic group it belongs to.} $\omega
_{V\otimes U}$ of $Sp(V\otimes U)$ which is convenient for us when we
consider the restriction of $\omega _{V\otimes U}$ to the subgroups $Sp$ and 
$O_{\beta }$. Let us identify the Lagrangian subspace $Y\otimes U$ of $%
V\otimes U$ with $Hom(X,U).$ This enables to realize (see Remark \ref%
{Schrodinger Model}) the representation $\omega _{V\otimes U}$ on the space
of functions 
\begin{equation*}
\mathcal{H=}L^{2}(Hom(X,U)).
\end{equation*}%
In this realization, the action of an element $A:Y\rightarrow X$ of the
Siegel unipotent $N$ of $Sp$ is given by 
\begin{equation}
\left( \omega _{V\otimes U}(A)f\right) (T)=\psi (\frac{1}{2}Tr(\beta
_{T}\circ A))f(T),  \label{Op-Model-N}
\end{equation}%
where for $T:X\rightarrow U$ we denote by $\beta _{T}:X\rightarrow X^{\ast
}=Y$ the quadratic form 
\begin{equation}
\beta _{T}(x,x^{\prime })=\beta (T(x),T(x^{\prime })),  \label{Betta-T}
\end{equation}%
and we denote by $Tr(\beta _{T}\circ A)$ the trace of the composite operator 
$Y\overset{A}{\rightarrow }X\overset{\beta _{T}}{\rightarrow }Y.$ In
addition, in this model the action of an element $r\in O_{\beta }$ is given
by 
\begin{equation}
\left( \omega _{V\otimes U}(r)f\right) (T)=\legendre{\det
(r)^{n}}{q}f(r^{-1}\circ T).  \label{Op-Model-O}
\end{equation}

\subsection{\textbf{The Eta Correspondence\label{EC}}}

Consider the oscillator representation $\omega _{V\otimes U}$ of $%
Sp(V\otimes U)$.

\begin{remark}
For the rest of this section we make the following choice. If $\psi $ is the
central character we use to define $\omega _{V\otimes U},$ then the
character $\psi _{\frac{1}{2}},$ $\psi _{\frac{1}{2}}(z)=\psi (\frac{1}{2}%
z), $ of $\mathbb{F}_{q}$ is the one we use in (\ref{SymXdual}) to identify $%
Sym^{2}(Y)$ and $\widehat{N}.$
\end{remark}

With this choice of parameters we can make the following precise statement:

\begin{proposition}
\label{Rank-as-rep-of-Sp}Assume that $\dim (U)=k<n.$ As a representation of $%
Sp,$ $\omega _{V\otimes U}$ is of rank $k$ and type\footnote{%
A rank $k$ form $B$ on $Y$ is of type $\mathcal{O}_{\beta }$ if $Y/rad(B)$
is isometric to $(U,\beta ).$} $\mathcal{O}_{\beta }.$
\end{proposition}

For a proof of Proposition \ref{Rank-as-rep-of-Sp}, see Appendix \ref%
{P-of-Proposition-Rank-as-rep-of-Sp}.\medskip

Now, consider the restriction, via the map (\ref{Embedding}), of $\omega
_{V\otimes U}$ to the product $Sp\times O_{\beta }.$ We decompose this
restriction into isotypic components for $O_{\beta }$:%
\begin{equation}
\omega _{V\otimes U|Sp\times O_{\beta }}\simeq \sum_{\tau \in Irr(O_{\beta
})}\Theta (\tau )\otimes \tau ,  \label{Res-DP}
\end{equation}%
where $\Theta (\tau )$ is a representation of $Sp.$ Although the factors $%
\Theta (\tau )$ in (\ref{Res-DP}) will in general not be irreducible, we can
say something about how they decompose. Let us denote by 
\begin{equation*}
Irr(Sp)_{k}\supset Irr(Sp)_{k\beta },
\end{equation*}%
the sets of (equivalence classes of) irreducible representations of $Sp$ of
rank $k,$ and of rank $k$ and type $\mathcal{O}_{\beta },$ respectively. The
next theorem---the main result of this note---announces that each $\Theta
(\tau )$ has a certain largest "chunk", which is in fact what we are
searching for.

\begin{theorem}[\textbf{Eta Correspondence}]
\label{EtaC}Assume that $\dim (U)=k<n.$ The following hold true:

\begin{enumerate}
\item \textbf{Rank }$k$\textbf{\ piece.} For each $\tau $ in $Irr(O_{\beta
}) $ the representation $\Theta (\tau )$ contains a unique irreducible
constituent $\eta (\tau )$ of rank $k$; all other constituents have rank
less than $k$.

\item \textbf{Injection. }The mapping $\tau \longmapsto \eta (\tau )$ gives
an embedding 
\begin{equation}
\eta :Irr(O_{\beta })\longrightarrow Irr(Sp)_{k\beta }.  \label{Eta}
\end{equation}

\item \label{Spectrum}\textbf{Spectrum. }The multiplicity of the orbit $%
\mathcal{O}_{\beta }$ in $\eta (\tau )_{|N}$ is $\dim (\tau ).$\smallskip
\end{enumerate}
\end{theorem}

For a proof of Theorem \ref{EtaC}, see Appendix \ref{P-of-Thm-EtaC}.\medskip

It also seems that this construction should produce all of the rank $k$
representations. We formulate this as a conjecture.

\begin{conjecture}[\textbf{Exhaustion}]
\label{Exhaustion}Assume that $\dim (U)=k<n.$ We have%
\begin{equation}
Irr(Sp)_{k}=\eta (Irr(O_{\beta +}))\text{ }\bigsqcup \text{ }\eta
(Irr(O_{\beta -})),  \label{Ex}
\end{equation}%
where $\beta +$ and $\beta -$ represent the two isomorphism classes of inner
products of rank $k$.
\end{conjecture}

\begin{remark}
Note that by (\ref{Eta}) the union in (\ref{Ex}) is indeed disjoint.
\end{remark}

Conjecture \ref{Exhaustion} is backed up by theoretical observations and
numerical computations---see Subsection \ref{Numerics} for illustration.

\begin{remark}[\textbf{The case} $\dim(U)=n$]
Proposition \ref{Rank-as-rep-of-Sp} and Theorem \ref{EtaC}, and their proofs, hold also in the case $\dim(U)=n$. However, due to Conjecture \ref{Exhaustion} we decided to formulate them with $k<n$.  
\end{remark}

We give now several additional remarks that, in particular, will clarify the novelty of
our main result, and will also explain why we decided to call (\ref{Eta})
the \underline{eta correspondence}.

\begin{remark}
We would like to comment that

\begin{description}
\item[(a) \textbf{Eta vs. Theta correspondence over local fields}] \textit{%
Considering the groups }$Sp$ \textit{and }$O_{\beta }$ \textit{over a local
field, one can associate, in a similar fashion as above, to every
irreducible representation }$\tau $ \textit{of }$O_{\beta },$ \textit{a
representation }$\Theta (\tau )$ of $Sp.$ It will in general not be
irreducible and the question is, what component to select from it? One
option is to take the "minimal" piece of $\Theta (\tau ).$ Indeed, it turns
out that $\Theta (\tau )$ has a unique irreducible quotient $\theta (\tau ).$
The assignment 
\begin{equation*}
\tau \mapsto \theta (\tau )\text{ \ - \ \textit{"minimal" piece,}}
\end{equation*}%
is the famous theta correspondence, which has been studied by many authors 
\cite{Ge, GT, H3, H5, K, MVW, P, R, Wa} for its usefulness in the theory of
automorphic forms. A second option is to take the \textit{"maximal" }piece
of $\Theta (\tau ).$ Indeed, repeating in the local field case, verbatim,
the scheme we proposed above, we find that $\Theta (\tau )$ has a largest
chunk in the form of a unique irreducible sub-representation $\eta (\tau )$
of rank $k$, which will equal $\theta (\tau )$ exactly when $\Theta (\tau )$
is irreducible. The assignment 
\begin{equation*}
\tau \mapsto \eta (\tau )\text{ \ - \ "maximal" piece,}
\end{equation*}%
is our eta correspondence (\ref{Eta}). Application of the new correspondence
to representation theory of classical groups over local fields, will be a
subject for future publications.\smallskip

\item[(b) \textbf{Eta Correspondence over finite fields}] As noted by
several authors (see \cite{AKT, AM, AMR, H2, Sr}, and in particular \cite%
{AKT} where the case of unipotent representations was considered) the theta
correspondence is not defined over finite fields\footnote{%
In fact, the attempt \cite{H2} to develop a duality theory over finite
fields preceded the one over the local fields \cite{H3}.}. The eta
correspondence comes as the appropriate construction in this case. This is
also the reason we use a related, although different, notation for the
correspondence (\ref{Eta}).
\end{description}
\end{remark}

Finally, we would like to make the following remark on the generality of our
work.

\begin{remark}[\textbf{Generalized Eta Correspondence}]
The notion of rank for the group $Sp_{2n}$ over local fields was described
in \cite{H4}. The theory for general classical groups over local fields was
developed by Li in \cite{L}, and it was extended to all semi-simple
algebraic groups over local fields by Salmasian in \cite{S}. The development
of the eta correspondence (\ref{Eta}) for all finite classical groups will
be discussed in future publications. For expositional purposes, in this note
we describe only the case of the finite symplectic groups.
\end{remark}

\subsection{\textbf{Numerical Justification for the Exhaustion Conjecture 
\label{Numerics}}}

Conjecture \ref{Exhaustion} is backed up by numerical data collected for the
groups $Sp_{6}(\mathbb{F}_{q}),$ $q=3,5,7,9,11,13$; $Sp_{8}(\mathbb{F}_{q})$%
, $q=3,5,$ and $Sp_{10}(\mathbb{F}_{3}).$ 
\begin{figure}[h]\centering
\includegraphics[
natheight=4.0836in, natwidth=12.9722in, height=1.95in, width=5.95in]
{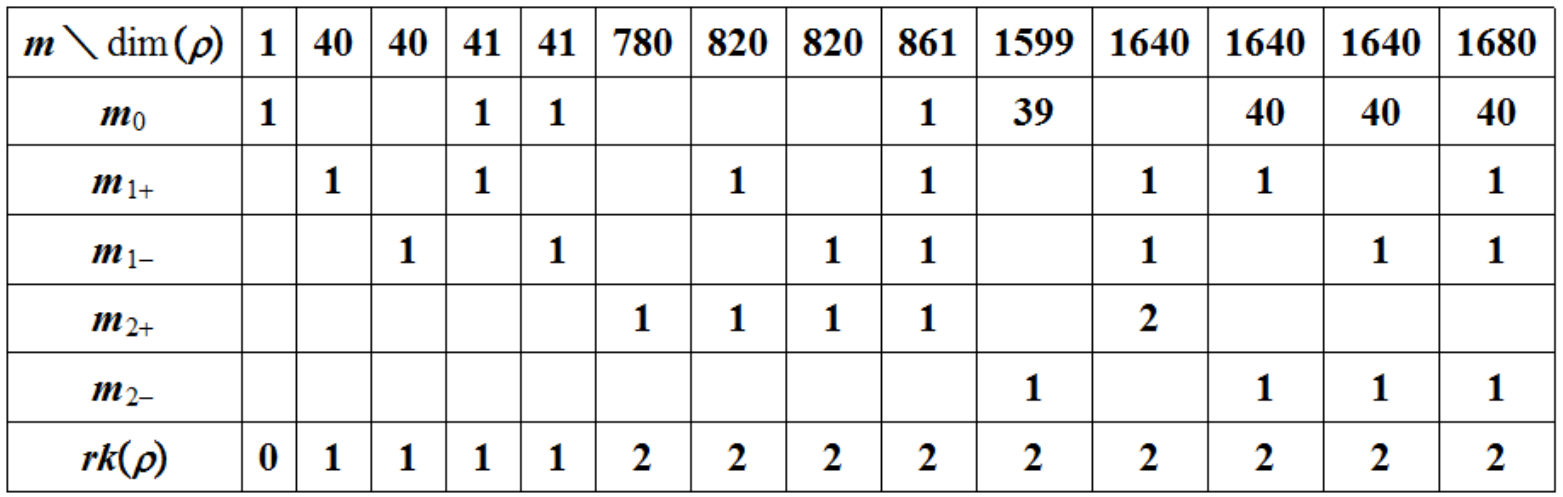}%
\caption{Multiplicities and rank for irreps of $Sp_{8}(\mathbb{F}_{3})$:
Ranks $k=0,$ $1,$ $2.$ }\label{Ranksp8_3-012}%
\end{figure}%
\begin{figure}[h]\centering
\includegraphics[
natheight=6.333in, natwidth=16.1668in, height=2.25in, width=5.9in]
{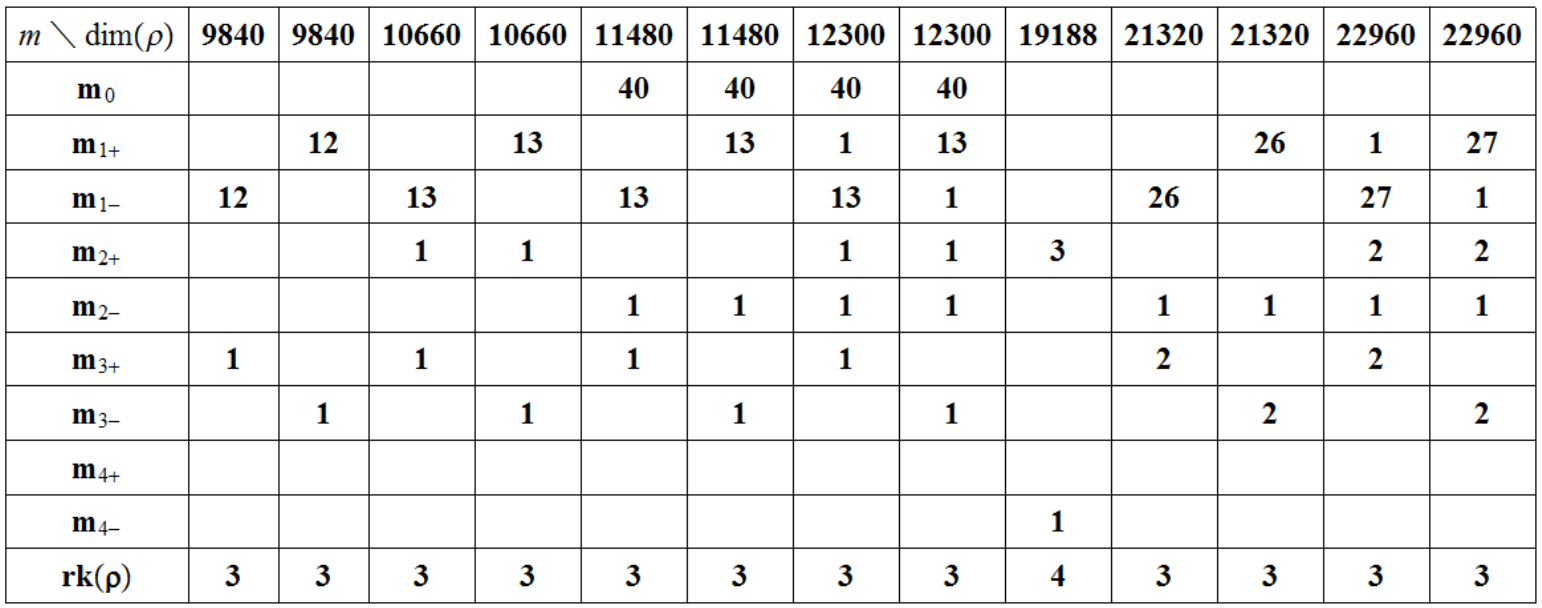}%
\caption{Multiplicities and rank for irreps of $Sp_{8}(\mathbb{F}_{3})$:
Ranks $k=3,$ $4.$ }\label{ranksp8_3-34a}%
\end{figure}%
\begin{figure}[h]\centering
\includegraphics[
natheight=6.333in, natwidth=16.0414in, height=2.25in, width=5.95in]
{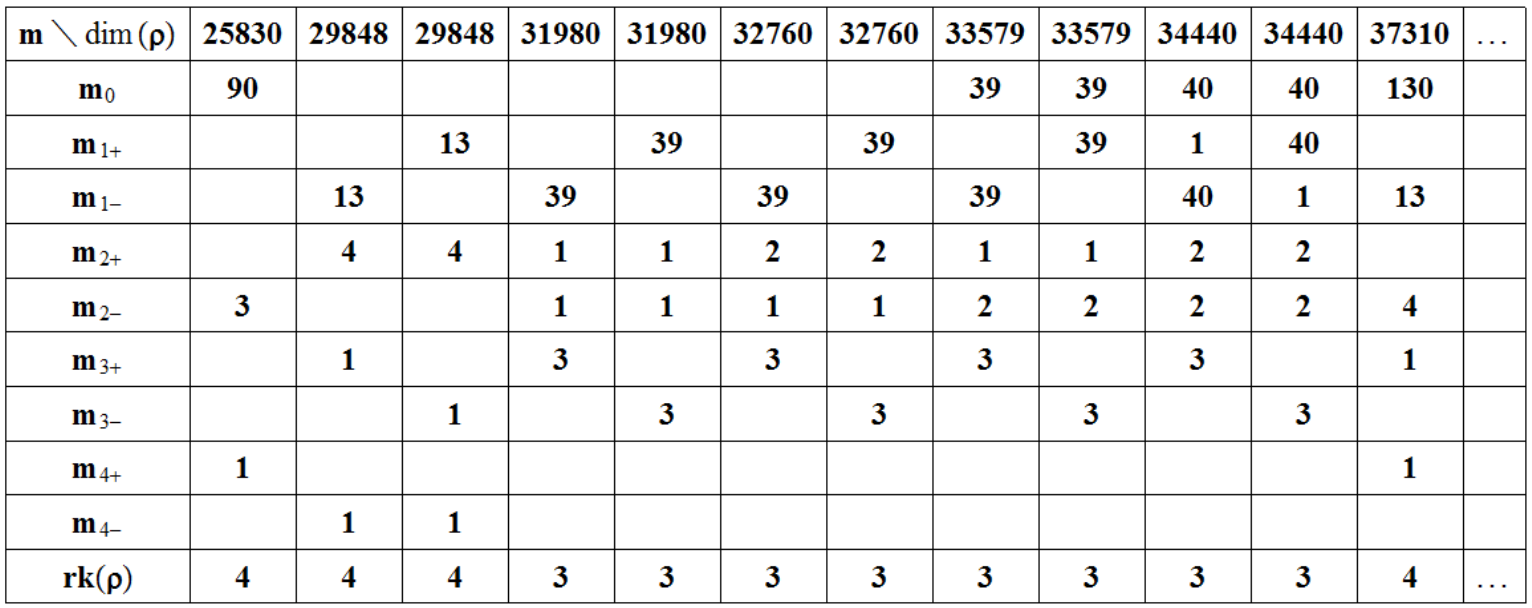}%
\caption{Multiplicities and rank for irreps of $Sp_{8}(\mathbb{F}_{3})$:
Ranks $k=3,$ $4.$ }\label{ranksp8_3-34b}%
\end{figure}%
Indeed, the Magma computations done with Cannon and Goldstein, for the
various sizes of symplectic groups, repeatedly confirm the assertion made in
the exhaustion conjecture, i.e., Identity (\ref{Ex}). For example, the
number of conjugacy classes in the orthogonal groups $O_{1+}(\mathbb{F}_{q})$
and $O_{1-}(\mathbb{F}_{q})$ together is $4,$ each contributes $2$ classes;
the number of conjugacy classes in the groups $O_{2+}(\mathbb{F}_{q})$ and $%
O_{2-}(\mathbb{F}_{q})$ together is $q+6,$ one contributes $\frac{q+5}{2}$
classes and the other $\frac{q+7}{2}$ classes. In addition, the number of
conjugacy classes in the groups $O_{3+}(\mathbb{F}_{q})$ and $O_{3-}(\mathbb{%
F}_{q})$ together is $4(q+2)$, each contributes $2(q+2)$ classes. Hence, the
computations of the multiplicities and rank for the groups $Sp_{6}(\mathbb{F}%
_{5})$ and $Sp_{8}(\mathbb{F}_{3})$ presented in Figures \ref{RankSp6_5} and %
\ref{Ranksp8_3-012}--\ref{ranksp8_3-34a}--\ref{ranksp8_3-34b}, respectively,
give the required numerical confirmation of (\ref{Ex}) in these cases.

\section{\textbf{Dimension of Rank }$k$\textbf{\ Representations}}

We would like to clarify the strong relationship between the dimension of a
representation of $Sp$ and its rank.

\subsection{\textbf{Dimension}}

\begin{figure}[h]\centering
\includegraphics[
natheight=4.7919in, natwidth=16.4306in, height=1.9441in, width=6.6011in]
{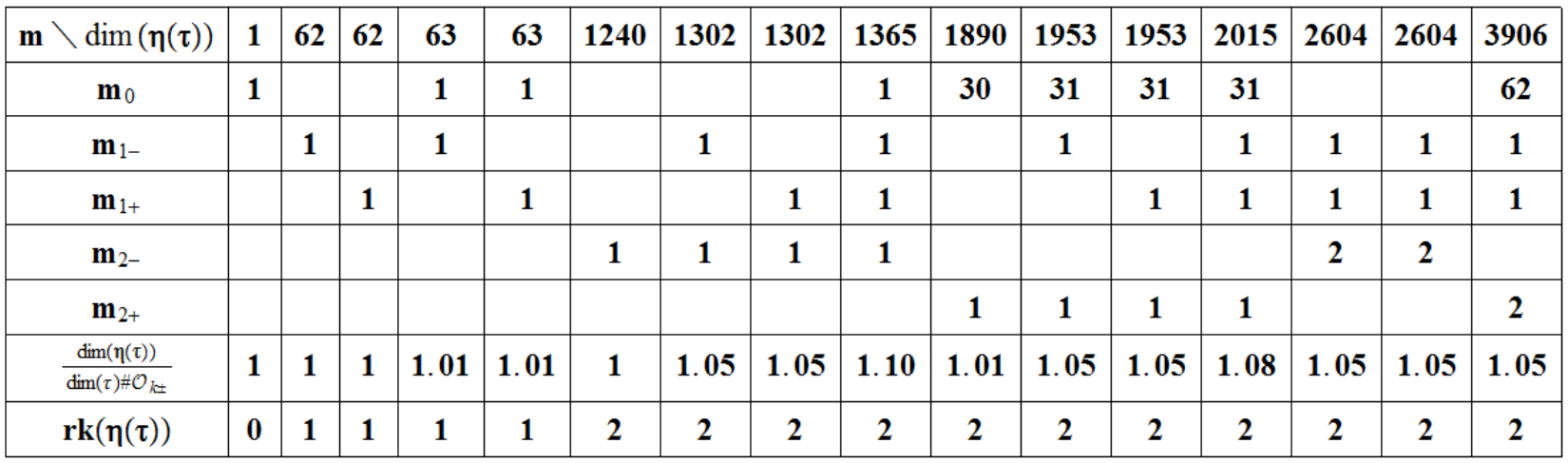}%
\caption{The relation between $\dim (\protect\tau )\cdot \#\mathcal{O}_{k\pm
}$ and $\dim (\protect\eta (\protect\tau ))$ for rank $k<3$ irreps of $%
Sp_{6}(\mathbb{F}_{5}).$}\label{dim-ranksp6_5}%
\end{figure}

Fix $k<n$ and consider a rank $k$ irreducible representation $\rho \in
Irr(Sp)_{k}.$ Let us assume that $\rho $ appears in the image of the eta
correspondence (\ref{Eta}). Namely, there exist $\tau \in Irr(O_{k\pm })$
such that $\rho =\eta (\tau )$---see Section \ref{EC}. Using Part \ref%
{Spectrum} of Theorem \ref{EtaC}, and the dimension formula (\ref{DimForm}),
we have 
\begin{equation}
\dim (\eta (\tau ))=\dim (\tau )\cdot \#\mathcal{O}_{k\pm
}+\sum_{r<k}\sum_{\pm }m_{r\pm }\cdot \#\mathcal{O}_{r\pm }.
\label{Dim-eta-tau-Form}
\end{equation}%
The point now is---see Figure \ref{dim-ranksp6_5} for illustration---that
the term $\dim (\tau )\cdot \#\mathcal{O}_{k\pm }$ dominates the right hand
side of (\ref{Dim-eta-tau-Form}). Indeed, we have

\begin{theorem}[\textbf{Dimension Estimate}]
\label{Dim}Let $\eta (\tau )$ be a rank $k<n$ irreducible representation of $%
Sp$ associated to an irreducible representation $\tau $ of $O_{k\pm }.$ Then 
\begin{equation}
1\leq \frac{\dim (\eta (\tau ))}{\dim (\tau )\#\mathcal{O}_{k\pm }}\leq 1+%
\frac{2+\varepsilon (q)}{q^{n-k+1}},\text{ \ with \ }\varepsilon (q)=O(1/q).
\label{DimEst}
\end{equation}
\end{theorem}

A proof of Theorem \ref{Dim} will be given in a sequel paper.\smallskip\ 

\begin{remark}
The term $\varepsilon (q)$ can be estimated explicitly. For example we have $%
\varepsilon (q)<2/q+4/q^{2}.$
\end{remark}

Theorem \ref{Dim} seems to substantially extend the current knowledge \cite%
{LS, N, Sh, TZ} on the dimensions of representations of the finite
symplectic groups (See Lemma 2.3. in \cite{M}).

\subsection{\textbf{Compatibility of Dimension and Rank}}

Although dimension tends to increase with rank, because of the factor $\dim
(\tau )$ in (\ref{Dim-eta-tau-Form}), it may happen, see Figures \ref%
{ranksp8_3-34a}--\ref{ranksp8_3-34b}, that a representation of rank $k$ has
larger dimension than one of rank $k+1$. However, for a given $k$, if $n$ is
large enough then the representations of rank $k$ will have smaller
dimension than those of rank $k+1$. For example, it seems that if $q$ is
sufficiently large, then for $k=1$, one can take $n=2$, and for $k=2,$ one
can take $n=3$---see Figures \ref{dim-ranksp6_5} and \ref{Ranksp8_3-012} for
illustration. In general, using Theorem \ref{Dim} and the known estimates on
the dimensions of the largest irreducible representations of the orthogonal
groups, we have the following result:

\begin{proposition}[\textbf{Compatibility of Dimension and Rank}]
\label{Compatibility}For sufficiently large $q,$ in the regime 
\begin{equation*}
k<2\sqrt{n}-1,
\end{equation*}%
the rank $k$ representations appearing in the image of the eta
correspondence (\ref{Eta}) always have smaller dimension then those of rank $%
k+1$.
\end{proposition}

The exact computation leading to a verification of Proposition \ref%
{Compatibility} will be given in a sequel paper.

\appendix

\section{\textbf{Proofs}}

\subsection{\textbf{Proof of Lemma} \protect\ref{Trivial}\label%
{P-of-Lemma-Trivial}}

\begin{proof}
If $Supp_{N}(\rho )=0$ then $\rho _{|N}$ is trivial. The Lemma now follows
from the well known fact that the $N$ conjugates generate the group $Sp$ 
\cite{A}.
\end{proof}

\subsection{\textbf{Proof of Proposition }\protect\ref{Oscillators}\label%
{P-of-Proposition-Oscillators}}

\begin{proof}
Consider the automorphism $\alpha _{s}:H\rightarrow H$ given by $\alpha
_{s}(v,z)=(v,s^{2}z).$ It can be extended to an automorphism of the
semi-direct product of $Sp$ with $H$, by letting it act trivially on $Sp$.
The equivalence of oscillator representations follows. The fact that for a \
non-square $\varepsilon \in \mathbb{F}_{q}^{\ast }$, the representations $%
\omega _{\psi }$ and $\omega _{\psi _{\varepsilon }}$ are not isomorphic,
can be verified using the realization given in Remark \ref{Schrodinger Model}%
. This completes the proof of the proposition.
\end{proof}

\subsection{\textbf{Proof of Proposition} \protect\ref{Rank-as-rep-of-Sp} 
\label{P-of-Proposition-Rank-as-rep-of-Sp}}

\begin{proof}
The proposition follows immediately from Equation (\ref{Op-Model-N}) in
Section \ref{Schrodinger-VU}.
\end{proof}

\section{\textbf{Proof of the Eta Correspondence Theorem}\label%
{P-of-Thm-EtaC}}

We give a proof of Theorem \ref{EtaC}\textbf{\ }that is an elementary
application of the \textit{double commutant theorem }\cite{We}.

\subsection{\textbf{The Double Commutant Theorem}}

We will use the following version:

\begin{theorem}[\textbf{Double Commutant Theorem}]
Let $W$ be a finite dimensional vector space. Let $\mathcal{A},\mathcal{%
\mathcal{A}}^{\prime }$ $\subset End(W)$ be two sub-algebras, such that

\begin{enumerate}
\item The algebra $\mathcal{A}$ acts semi-simply on $W.$

\item Each of $\mathcal{A}$ and $\mathcal{\mathcal{A}}^{\prime }$ is the
full commutant of the other in $End(W).$
\end{enumerate}

Then $\mathcal{A}^{\prime }$ acts semi-simply on $W,$ and as a
representation of\ $\mathcal{A\otimes \mathcal{A}}^{\prime }$ we have%
\begin{equation*}
W=\bigoplus\limits_{i\in I}W_{i}\otimes W_{i}^{\prime },
\end{equation*}%
where $W_{i}$ are all the irreducible representations of $\mathcal{A}$, and $%
W_{i}^{\prime }$ are all the irreducible representations of $\mathcal{A}%
^{\prime }.$ In particular, we have a bijection between irreducible
representations of $\mathcal{A}$ and $\mathcal{A}^{\prime },$ and moreover,
every isotypic component for $\mathcal{A}$ is an irreducible representation
of $\mathcal{A\otimes \mathcal{A}}^{\prime }.$
\end{theorem}

\subsection{\textbf{Preliminaries}}

Let us start with several preliminary steps. We work with the Schr\"{o}%
dinger model of $\omega _{V\otimes U}$ appearing in Section \ref%
{Schrodinger-VU}. It is realized on the space 
\begin{equation}
\mathcal{H}=L^{2}(Hom(X,U)),  \label{OsRep}
\end{equation}%
and there, the actions of an element $A$ of the Siegel unipotent radical $N$ 
$\subset Sp,$ and an element $r\in O_{\beta },$ are given by Formulas (\ref%
{Op-Model-N}) and (\ref{Op-Model-O}), respectively. In particular, we have

\begin{claim}
Every character appearing in the restriction of $\omega _{V\otimes U}$ to $N$
is of the form $\psi _{\beta _{T}}$ for some $T\in Hom(X,U).$ Moreover, we
have $rank(\beta _{T})=k$ iff $T$ is onto.
\end{claim}

For the rest of the section, we fix a transformation $T:X\twoheadrightarrow
U $ which is \underline{onto} and consider the character subspace 
\begin{equation*}
\mathcal{H}^{\psi _{\beta _{T}}}=\{f\in \mathcal{H};\text{ }\omega
_{V\otimes U}(A)f\text{ }=\psi _{\beta _{T}}(A)f\text{, \ }A\in N\}.
\end{equation*}%
We would like to have a better description of the space $\mathcal{H}^{\psi
_{\beta _{T}}}$. The orthogonal group $O_{\beta }$ acts naturally on $%
Hom(X,U)$ and we denote by $\mathcal{O}_{T}$ the orbit of $T$ under this
action.

\begin{proposition}
\label{Regular}We have $\mathcal{H}^{\psi _{\beta _{T}}}=L^{2}(\mathcal{O}%
_{T})$ the space of functions on $\mathcal{O}_{T}.$
\end{proposition}

For a proof of Proposition \ref{Regular}, see Section \ref%
{P-Proposition-Regular}.\medskip

Note that, because $T$ is onto, the action of $O_{\beta }$ on $\mathcal{O}%
_{T}$ is free. In particular, we can identify $\mathcal{O}_{T}$ with $%
O_{\beta },$ and the Peter--Weyl theorem \cite{Se} for the regular
representation implies

\begin{corollary}
Under the action of $O_{\beta },$ the space $\mathcal{H}^{\psi _{\beta
_{T}}} $ decomposes as 
\begin{equation}
\mathcal{H}^{\psi _{\beta _{T}}}\simeq \bigoplus\limits_{\tau \in
Irr(O_{\beta })}\dim (\tau )\tau .  \label{Reg-of-Obeta}
\end{equation}
\end{corollary}

We would like now to describe the commutant of $O_{\beta }$ in $End(\mathcal{%
H}^{\psi _{\beta _{T}}}).$ Considering the group 
\begin{equation*}
G_{\beta _{T}}=Stab_{GL(X)}(\beta _{T}),
\end{equation*}
of automorphisms of $X$ that stabilize the form $\beta _{T},$ we obtain two
commuting actions%
\begin{equation*}
O_{\beta }\curvearrowright \mathcal{H}^{\psi _{\beta _{T}}}\curvearrowleft
G_{\beta _{T}}.
\end{equation*}%
Moreover, we have

\begin{proposition}
\label{Commutant}\textit{The groups }$O_{\beta }$\textit{\ and }$G_{\beta
_{T}}$\textit{\ generate each other's commutant in }$End(\mathcal{H}^{\psi
_{\beta _{T}}}).$
\end{proposition}

For a proof of Proposition \ref{Commutant}, see Section \ref%
{P-Proposition-Commutant}.\smallskip

\subsection{\textbf{Proof of Theorem \protect\ref{EtaC}}}

\begin{proof}
Write 
\begin{equation*}
\Theta (\tau )\simeq \sum \eta _{i}(\tau ),
\end{equation*}%
for various irreducible representations $\eta _{i}(\tau )$ of $Sp$. Then%
\begin{equation}
\Theta (\tau )^{\psi _{\beta _{T}}}\simeq \sum \eta _{i}(\tau )^{\psi
_{\beta _{T}}}.  \label{D}
\end{equation}%
In addition, $\Theta (\tau )^{\psi _{\beta _{T}}}$ is a $G_{\beta _{T}}$%
-module, and so is each\ $\eta _{i}(\tau )^{\psi _{\beta _{T}}}.$ Hence,
Identity (\ref{D}) gives a decomposition of $\Theta (\tau )^{\psi _{\beta
_{T}}}$ into (not necessarily irreducible) submodules for $G_{\beta _{T}}$.
But Proposition \ref{Commutant} together with the Double Commutant Theorem
says that $\Theta (\tau )^{\psi _{\beta _{T}}}$ is irreducible as a $%
G_{\beta _{T}}$-module. Therefore, exactly one of the $\eta _{i}(\tau
)^{\psi _{\beta _{T}}}$ will be non-zero, and it defines an irreducible
representation of $G_{\beta _{T}}$, which has dimension equal to $\dim (\tau
)$, by equation (\ref{Reg-of-Obeta}). To conclude, there exists a unique
irreducible sub-representation $\eta (\tau )<$ $\Theta (\tau )$ of rank $k$
and type $\mathcal{O}_{\beta _{T}},$ the multiplicity of the orbit $\mathcal{%
O}_{\beta }$ in $\eta (\tau )_{|N}$ is $\dim (\tau ),$ and finally, the
Double Commutant Theorem implies that for $\tau \ncong \tau ^{\prime }$ in $%
Irr(O_{\beta }),$ we have $\eta (\tau )\ncong \eta (\tau ^{\prime }).$ This
completes the proof of Theorem \ref{EtaC}.
\end{proof}

\subsection{\textbf{Proofs}}

\subsubsection{\textbf{Proof of Proposition \protect\ref{Regular}}\label%
{P-Proposition-Regular}}

\begin{proof}
Using the delta basis $\{\delta _{T};$ $T\in Hom(X,U)\},$ we can verify
Claim \ref{Regular}, by showing that if $\beta _{T^{\prime }}=\beta _{T}$
then there exists $r\in O_{\beta }$ such that $T^{\prime }=r\circ T.$
Indeed, let $r$ be the composition 
\begin{equation*}
U\widetilde{\longrightarrow }X/rad(\beta _{T})\widetilde{\longrightarrow }U,
\end{equation*}%
where the first and second isomorphisms are these induced by $T^{\prime },$
and $T,$ respectively, and $rad(\beta _{T})$ is the radical of $\beta _{T}.$
This completes the proof of Proposition \ref{Regular}.
\end{proof}

\subsubsection{\textbf{Proof of Proposition \protect\ref{Commutant}\label%
{P-Proposition-Commutant}}}

\begin{proof}
To verify this assertion, note that we have a short exact sequence%
\begin{equation}
1\rightarrow N_{k,n-k}\rightarrow G_{\beta _{T}}\rightarrow O(X/rad(\beta
_{T}))\times GL(Z)\rightarrow 1,  \label{SES}
\end{equation}%
where $Z=\ker (T)=rad(\beta _{T})$, $O(X/rad(\beta _{T}))$ is the orthogonal
group of $X/rad(\beta _{T}),$ and $N_{k,n-k}$ is the appropriate unipotent
group. The group $O(X/rad(\beta _{T}))$ acts simply transitively on the
orbit $\mathcal{O}_{T}$, as does the group $O_{\beta }$, and these two
actions commute with each other. If we use the map $r\mapsto r^{-1}\circ T$
to identify $O_{\beta }$ with $\mathcal{O}_{T}$, then the action of $%
O_{\beta }$ becomes the action of $O_{\beta }$ on itself by left
translation, and the action of $O(X/rad(\beta _{T}))$ can be identified with
the action of $O_{\beta }$ on itself by right translation. By the Peter-Weyl
Theorem \cite{Se}, we conclude that the groups $O(X/rad(\beta _{T}))$ and $%
O_{\beta }$ generate mutual commutants in the operators on $L^{2}(\mathcal{O}%
_{T})\simeq \mathcal{H}^{\psi _{\beta _{T}}}$. \underline{A fortiori} the
groups $G_{\beta _{T}}$ and $O_{\beta }$ generate mutual commutants on $%
L^{2}(\mathcal{O}_{T})$. This completes the proof of Proposition \ref%
{Commutant}.
\end{proof}

\end{document}